\documentclass{commat}

\usepackage{graphicx}

\title{%
    Dynamic behavior of a railway track under a moving wheel load modelled as a sinusoidal pulse
    }

\author{%
    Nouzha Lamdouar, Chakir Tajani and Mohammed Touati
    }

\affiliation{
    \address{Nouzha Lamdouar --
    Civil Department, Mohammadia School of Engineers, Mohammed V University, Morocco
        }
    \email{%
    nlamdouar@gmail.com
    }
    \address{Chakir Tajani --
    Department of Mathematics, Polydisciplinary Faculty of Larache, Abdelmalek Essaadi University, Morocco
        }
    \email{%
    chakir_tajani@hotmail.fr
    }
    \address{Mohammed Touati --
    Civil Department, Mohammadia School of Engineers, Mohammed V University, Morocco
        }
    \email{%
    mohammed.touati87@gmail.com
    }
    }

\abstract{%
    The aim of this paper is to evaluate the train/track induced loads on the substructure by modelling the wheel, at each instant, as a moving sinusoidal pulse applied in a very short period of time. This assumption has the advantage of being more realistic as it reduces the impact of time on the load definition.
    To that end, mass, stiffness, and dumping matrices of an elementary section of track will be determined. As a result, the equations of motion of a section of track subjected to a sinusoidal pulse and a rectangular pulse respectively is concluded. Two numerical methods of resolution of that equation, depending on the nature of the dumping matrix, will be presented. The computation results will be compared in order to conclude about the relevance of that load model. This approach is used in order to assess the nature and the value of the loads received by the substructure.
    }



\VOLUME{32}
\NUMBER{1}
\YEAR{2024}
\firstpage{35}
\DOI{https://doi.org/10.46298/cm.10774}

\begin{paper}

\section{Introduction}
Various theoretical and experimental researches have been performed in order to assess train/track induced loads on the substructure. Mohammed Touati and al. \cite{Key1} determined the loads induced by a non-linear 3D multi-body modelled train on the track with taking into account wheel/rail contact properties and track irregularities. Yang Xinwen and al. \cite{Key2} concluded, through a vehicle-track-subgrade coupling dynamic theory and finite element method, about the train/track induced loads on each layer of the substructure. As an experimental study, Al Shaer and al. \cite{Key3} presented the dynamic behavior of a portion of ballasted railway track subjected to cyclic loads in substitution of a moving wheelset. In conclusion, the dynamics behavior of the substructure is widely studied in the literature (\cite{Key4}, \cite{Key5}, \cite{Key6}, \cite{Key7}, \cite{Key8}) based on the train/track coupling model.\\
Actually, even if modelling a wheel load as a rectangular pulse is a common assumption, real measurements don't show the same shape. In fact, ONCF (Moroccan railway network manager) has many tools that record wheel pulse like GOTCHA. This system shows that the shape of the load has never been rectangular, but it's more likely compared to a sinusoidal pulse. Then, this paper deals with evaluating train/track induced loads on the substructure by proposing a new approach when it comes to modelling the shape of the wheel impact. Indeed, it's common to consider a moving load as a rectangular impulse applied on the nodes of a mesh structure in each period of time depending on signal sampling. This paper shows that assuming the wheel load as a sinusoidal pulse may reduce the impact of the period of time of its application and, consequently, minimize the loads induced on the substructure oversized by the common assumption. In that matter, a finite element model of the track will be presented and the numerical results will be compared.

\section{Track elementary section modeling} 
\subsection{Determination of mass, stiffness et dumping matrices} 
Let’s assume a portion of ballasted track composed of two elements of rail considered as a continuous Euler-Bernoulli beam, fixed to two sleepers by a couples of springs/dampers representing the railpads. The ballast is modelled as a couples of springs/dampers under each sleeper (Figure 1).
\begin{figure}
	\begin{center}	
		\includegraphics[width=6.8cm,height=6.cm]{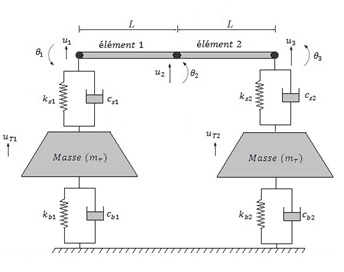}
	\end{center}
	\caption{Elementary track modelling}
\end{figure}

The displacement vector is written as:
\[
U=\left[u_{1}, \theta_{1}, u_{2}, \theta_{2}, u_{3}, \theta_{3}, u_{T 1}, u_{T 2}\right]
\]
The effective mass and the stiffness matrices of an element of rail \cite{Key9}, are given by:
\[M_{r}=\left(\rho_{r} A_{r} L / 420\right)\left[\begin{array}{cccccc}
	156 & 22 L & 54 & -13 L & 0 & 0 \\
	22 L & 4 L^{2} & 13 L & -3 L^{2} & 0 & 0 \\
	54 & 13 L & 312 & 0 & 54 & -13 L \\
	-13 L & -3 L^{2} & 0 & 8 L^{2} & 13 L & -3 L^{2} \\
	0 & 0 & 54 & 13 L & 156 & -22 L \\
	0 & 0 & -13 L & -3 L^{2} & -22 L & 4 L^{2}
\end{array}\right] \]
\[K_{r}=\left(E_{r} I_{r} / L^{3}\right)\left[\begin{array}{cccccc}
	12 & 6 L & -12 & 6 L & 0 & 0 \\
	6 L & 4 L^{2} & -6 L & 2 L^{2} & 0 & 0 \\
	-12 & -6 L & 24 & 0 & -12 & 6 L \\
	6 L & 2 L^{2} & 0 & 8 L^{2} & -6 L & 2 L^{2} \\
	0 & 0 & -12 & -6 L & 12 & -6 L \\
	0 & 0 & 6 L & 2 L^{2} & -6 L & 4 L^{2}
\end{array}\right]\]
where $\rho_{r}$ is the density of the rail, $A_{r}$ is the surface of the rail section, $E_{r}$ is Young
modulus, and $I_{r}$ is the rail moment of inertia. The dumping matrix of the rail is obtained as a linear combination of mass and stiffness matrices by assuming that the displacements $u_{1}$ and $u_{3}$ are completely dumped by the effect of railpads. \\
Therefore, the dumping matrix is written as:
$$
C_{r}^{*}=a_{0} \cdot M_{r}^{*}+a_{1} \cdot K_{r}^{*}
$$
where,
\[
M_{r}^{*}=\left(\rho_{r} A_{r} L / 420\right)\left[\begin{array}{cccc}
	4 L^{2} & 13 L & -3 L^{2} & 0 \\
	13 L & 312 & 0 & -13 L \\
	-3 L^{2} & 0 & 8 L^{2} & -3 L^{2} \\
	0 & 13 L & -3 L^{2} & 4 L^{2}
\end{array}\right]
\]
\[K_{r}^{*}=\left(E_{r} I_{r} / L^{3}\right)\left[\begin{array}{cccc}4 L^{2} & -6 L & 2 L^{2} & 0 \\ -6 L & 24 & 0 & 6 L \\ 2 L^{2} & 0 & 8 L^{2} & 2 L^{2} \\ 0 & 6 L & 2 L^{2} & 4 L^{2}\end{array}\right]\]
$a_0$ and $a_{1}$ are concluded from the equation:
$$\left[\begin{array}{l}a_{0} \\ a_{1}\end{array}\right]=\left(2 \omega_{1} \omega_{2} /\left(\omega_{2}^{2}-\omega_{1}^{2}\right)\right)\left[\begin{array}{cc}\omega_{2} & -\omega_{1} \\ -1 / \omega_{2} & 1 / \omega_{1}\end{array}\right]\left[\begin{array}{l}\zeta_{1} \\ \zeta_{2}\end{array}\right]$$
where $\omega_{i}{ }^{2}$, $(i=1,2)$ are the eigenvalues associated to the vibration of the rail described by the matrices $M_{r}{ }^{*}$ and $K_{r}{ }^{*},$ and $\zeta_{i}$, $(i=1,2)$ are the dumping ratios according to the first and second modes. 

In one hand, the equation of motion of the rail is written as:
\begin{equation}
	M_{r} \ddot{U}^{*}+C_{r} \dot{U}^{*}+K_{r} U^{*}=F
\end{equation}
where $C_{r}$ is the transformation of the matrix $C_{r}^{*}$ in the base $U^{*}$, and $U^{*}$ is defined by:
$$
U^{*}=\left[u_{1}, \theta_{1}, u_{2}, \theta_{2}, u_{3}, \theta_{3}\right]
$$
$F$ is given by:
$$
F=\left[\begin{array}{c}
	-k_{s}\left(u_{1}-u_{T 1}\right)-c_{s}\left(\dot{u}_{1}-\dot{u}_{T 1}\right) \\
	0 \\
	0 \\
	0 \\
	-k_{s}\left(u_{3}-u_{T 3}\right)-c_{s}\left(\dot{u}_{3}-\dot{u}_{T 3}\right) \\
	0
\end{array}\right]
$$
In the other hand, the equations of motion of the sleepers are written as:
\begin{equation}
	\left\{\begin{array}{c}
		m_{T} \ddot{u}_{T 1}=k_{s}\left(u_{1}-u_{T 1}\right)+c_{s}\left(\dot{u}_{1}-\dot{u}_{T 1}\right)-k_{b} u_{T 1}-c_{b} \dot{u}_{T 1} \\
		m_{T} \ddot{u}_{T 2}=k_{s}\left(u_{3}-u_{T 2}\right)+c_{s}\left(\dot{u}_{3}-\dot{u}_{T 2}\right)-k_{b} u_{T 2}-c_{b} \dot{u}_{T 2}
	\end{array}\right.
\end{equation}

From (1) and (2), we may conclude about the equation of motion of the track elementary section as it's modelled. It's written as:
$$
M \ddot{U}+C \dot{U}+K U=0
$$
where $M$, $C$ and $K$ are the mass, dumping, and the stiffness of the track elementary section respectively.

\subsection{Numerical application}
Let's assume a track elementary section characterized by the data given in table 1 (we can refer to (\cite{Key10}, \cite{Key11}, \cite{Key12}).

\begin{table}
	\begin{center}
		\begin{tabular}{|p{1.0in}|p{2.4in}|p{1.in}|}
			 \hline 	Symbol & Quantity & Value \\
			
			 \hline $\rho$${}_{r}$ & Rail density (kg/m${}^{3}$) & 7850 \\ 
			 
			\hline  \textit{A${}_{r}$} & Rail section surface (cm²) & 76.70 \\ 
			
			\hline \textit{E${}_{r}$} & Young modulus of the rail (GPa) & 210 \\
			
			 \hline  \textit{I${}_{r}$} & Rail moment of inertia (cm${}^{4}$) & 3038.6 \\
			  
			\hline 	\textit{m${}_{T}$} & Sleeper mass (kg) & 90.84 \\
			
			 \hline k${}_{s}$\textit{} & Railpad stiffness (MN/m) & 90 \\ 
			
			\hline 	$c{}_{s}$ & Railpad damping (kN.s/m) & 30 \\
			
			 \hline \textit{k${}_{b}$} & Ballast stiffness (MN/m) & 25.5 \\ 
			
			\hline \textit{c${}_{b}$} & Ballast damping (kN.s/m) & 40 \\
			
			 \hline $\zeta$\textit{} & Rail dumping ratio & 5\% \\ 
			
			\hline 
		\end{tabular}
		\caption{Track properties}
	\end{center}
\end{table}

The figure 2 illustrates the evolution of natural frequencies according to vibration modes. It shows that:

\begin{itemize}
	\item The frequencies of the $1^{st}$ and $2^{nd}$ modes correspond to a movement in phase between rail and sleepers. It's equal to 81.62 Hz;
	\item  The frequency of the $3^{rd}$ mode corresponds to a movement in opposition of phase between rail and sleepers. It's equal to 381.1 Hz.
\end{itemize}
\begin{figure}
	\begin{center}	
		\includegraphics[scale=0.85]{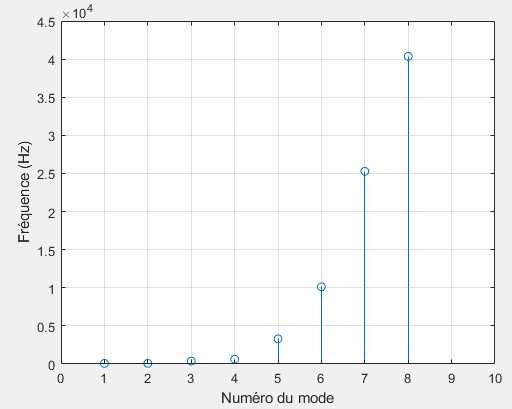} 
	\end{center}
	\caption{Natural frequencies of an elementary track section}
\end{figure}

\section{Track response to a rectangular and a sinusoidal pulses}
\subsection{Description of the studied track} 

Let's assume a section of track composed of $N$ track elementary sections subjected to an external load $F$ as it's shown in figure 3.

\begin{figure}
	\begin{center}	
		\includegraphics[scale=0.75]{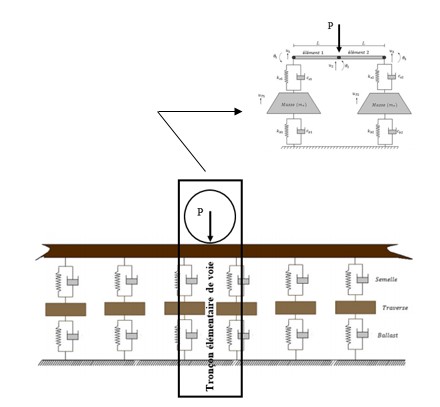} 	
	\end{center}
	\caption{  Track section modelling }	
\end{figure}

The number of degrees of freedom is given by:
$$N_{dof}=8N-3(N-1)$$
The displacement vector is written as:
\[U=\left[\begin{array}{c} {\vdots } \\ {u_{j,k} } \\ {\vdots } \end{array}\right]\] 
where,
\[\left\{\begin{array}{llll}
	u_{j,k} =u_{j,k}^{*} & where\quad k\in \left[1,8\right] & if & j=1 \\
	u_{j,k} =u_{j,k}^{*} & where\quad k\in \left[3,4,5,6,8\right] & if & j\ne 1 \end{array}\right. \] 
and,
\[U_{j}^{*} =\left[u_{j,1} ,\theta _{j,1} ,u_{j,2} ,\theta _{j,2} ,u_{j,3} ,\theta _{j,3} ,u_{j,T1} ,u_{j,T2} \right]\] 
$j$ refers to the element's number.

The mass, stiffness and dumping matrices in the base $U$ are obtained by assembling those of a track elementary section determined earlier.

The vector of loads is defined by:
\[F=\left[\begin{array}{c} {\vdots } \\ {f_{j} } \\ {\vdots } \end{array}\right]\] 
where,
\[\left\{\begin{array}{l}
	N \text { is even }
	\left\{\begin{array}{l}
		N=2 \left\{ \begin{array}{ll}
			f_{j}=P  & \text { if } \quad j=5 \\  
			f_{j}=0   & \text { else }  \end{array}\right. \\\\
		
		N \neq 2\left\{  \begin{array}{ll}
			f_{j}=P & \text { if } \quad j=(5 N / 2)+1 \\
			f_{j}=0 & \text { else } \end{array}\right. \\
	\end{array}\right. \\\\
	
	N \text { is uneven }\left\{\begin{array}{ll}
		f_{j}=P & \text { if } j=(5(N+1) / 2)-1 \\
		f_{j}=0  & \text { else }
	\end{array}\right.
\end{array}\right.\]
$P$ is a rectangular or a sinusoidal load given as:
\begin{itemize}
	\item  Sinusoidal pulse:
	\[\left\{\begin{array}{ll} 
		P=P_{0} \sin \omega t & if \quad t\le t_{d} \\ 
		P=0 & else \end{array}\right. \] 
	
	\item  Rectangular pulse:
	
	\[\left\{\begin{array}{ll} P=P_{0} & if\quad t\le t_{d}  \\
		P=0      & else \end{array}\right. \] 
\end{itemize}
Its shape is shown in the figure 4.
\begin{figure}
	\begin{center}	
		\includegraphics[scale=0.5]{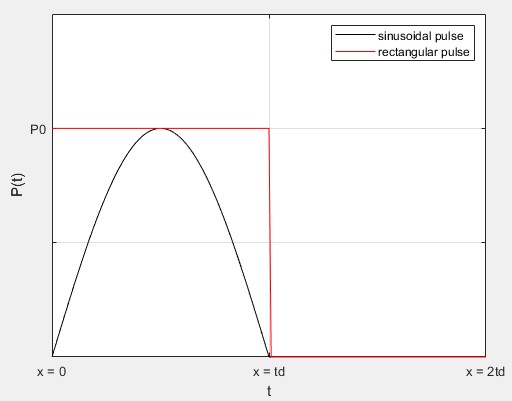} \\	
	\end{center}
	\caption{   Sinusoidal and rectangular pulses over a period of $t_d$ }	
\end{figure}

\subsection{ Description of the methods of resolution}

The dynamic behavior of the section of track may be analyzed by modal superposition if the dumping matrix verifies orthogonality properties. That method is used in particular for an undumped system. In that case, the equation of motion is reduced to:
\[M\ddot{U}+KU=F.\] 

Let's assume that ${\omega_{i}^{2}}$ are the eigenvalues associated to the track vibration. We note $\mathrm{\{}$$\mathrm{\phi}$${}_{i}$$\mathrm{\}}$ the normalized eigenvectors related to ${\omega_{i}^{2}}$. Therefore, the equation of motion is written as:
\begin{equation} 
	\ddot{Z}+diag(\omega _{i}^{2} )Z=\phi ^{T} F 
\end{equation} 
where diag(${\omega_{i}^{2}}$) is a diagonal matrix of the eigenvalues and:
\[U=\Phi .Z\] 
The system of equations (3) is uncoupled where each equation is written as:
\[\ddot{z}_{i} +\omega _{i}^{2} z_{i} =\Phi _{j,i} P(t)\] 
The resolution of that equation is given by DUHAMEL integral:
\[z_{i} (t)=(1/\omega _{i} )\int _{0}^{t}\Phi _{j,i} P(\tau )\sin \omega _{i} (t-\tau )d\tau  \] 
Therefore, the solution for a sinusoidal pulse load is given as:
\[z_{i}(t)=\left\{\begin{array}{lll}
	(\Phi_{j,i} P_{0} /\omega _{i}^{2} ).(1/(1-\beta ^{2} ))(\sin \omega t-\beta \sin \omega_{i} t) & if & t\le t_{d}  \\
	(\dot{z}_{i} (t_{d} )/\omega _{i} )\sin \omega _{i} (t-t_{d} )+z_{i} (t_{d} )\cos \omega _{i} (t-t_{d} ) & if & t\ge t_{d}  \end{array}\right. \]
where,
\[\beta =\omega /\omega _{i} \] 
and the solution for a rectangular pulse load is given as:
\[z_{i} (t)=\left\{\begin{array}{lll}
	(\Phi _{j,i} P_{0} /\omega _{i}^{2} )(1-\cos \omega _{i} t) & if\ & t\le t_{d}  \\
	(\Phi _{j,i} P_{0} /\omega _{i}^{2} )(\cos \omega _{i} (t-t_{d} )-\cos \omega _{i} t & if  & t\ge t_{d}  \end{array}\right. \] 

The figure 5 shows the response $z(t)$ to a sinusoidal and a rectangular pulse. It's obvious that in the forced phase, the maximum rectangular response is higher than the maximum sinusoidal response.

\begin{figure}
	\begin{center}	
		\includegraphics[width=7.2cm,height=6.cm]{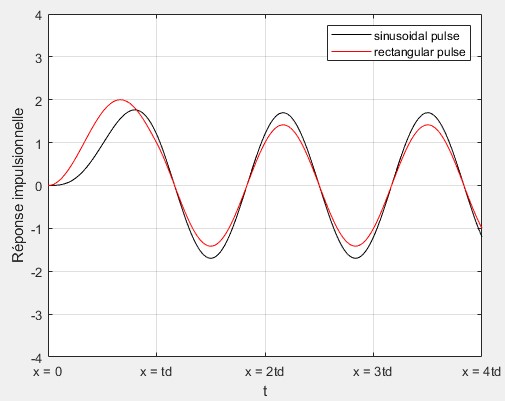}
	\end{center}
	\caption{$z(t)$ response to a rectangular and sinusoidal pulse $t_{d}/T = 0.75$}	
\end{figure}

In general, the dumping matrix doesn't verify the orthogonality characteristics. Therefore, the modal superposition method is substituted by the following method.

The equation of motion can be written as:
\begin{equation} 
	\ddot{Z}+\phi ^{T} C\phi .\dot{Z}+diag(\omega_{i}^{2} ).Z=\phi ^{T} F 
\end{equation} 
where diag($\omega$${}_{i}$${}^{2}$) and $\mathrm{\phi}$ are defined earlier. Knowing that:
\begin{equation}  
	\dot{Z}-\dot{Z}=0 
\end{equation} 
(4) and (5) could be written as:
\begin{equation} 
	\dot{Y}=D.Y+F^{*}  
\end{equation} 
where,
\[Y=\left[\begin{array}{c} {Z} \\ {\dot{Z}} \end{array}\right],\quad D=A^{-1} B{\kern 1pt} ,\quad F^{*} =A^{-1} \left[\begin{array}{c} {\phi ^{T} F} \\ {0} \end{array}\right]\] 
and,
\[A=\left[\begin{array}{cc} {\phi ^{T} C\phi } & {I} \\ {I} & {0} \end{array}\right],\quad B=\left[\begin{array}{cc} {diag(\omega _{i}^{2} )} & {0} \\ {0} & {-I} \end{array}\right]\]

Let's assume that $\mathrm{\{}$$\lambda$${}_{i}$$\mathrm{\}}$ are the eigenvalues associated to the matrix $D$. We note $\mathrm{\{}$$\psi$${}_{i}$$\mathrm{\}}$ the normalized eigenvectors related to $\mathrm{\{}$$\omega$${}_{i}$${}^{2}$$\mathrm{\}}$. We define $X(t)$ as:
\[Z=\psi .X\] 
The equation (6) is written as:
\begin{equation} 
	\dot{X}=diag(\lambda_{i} ).X+\psi ^{-1} F^{*}  
\end{equation} 
The system of equations (7) is uncoupled where each equation is written as:
\begin{equation} 
	\dot{x}_{i} (t)=a_{i} .x_{i} (t)+b_{i} .P 
\end{equation} 
where,
\[a_{i} =\lambda _{i} \quad and\quad b_{i} =\chi _{i} \] 
and,
\[\chi =\psi ^{-1} \left[\begin{array}{cc} {\phi ^{T} } & {0} \\ {0} & {0} \end{array}\right]\] 

The resolution of the equation (8) gives:

\begin{itemize}
	\item  Sinusoidal pulse:
	
	\[x_{i} (t)=\left\{\begin{array}{lll}
		\frac{b_{i} P\omega }{a_{i}^{2} +\omega ^{2} } e^{a_{i} t} -\frac{a_{i} b_{i} P}{a_{i}^{2} +\omega ^{2} } \sin \omega t-\frac{b_{i} P\omega }{a_{i}^{2} +\omega ^{2} } \cos \omega t & if  & t\le t_{d}  \\
		x_{i} (t_{d} )e^{a_{i} (t-t_{d} )} & if & t\ge t_{d}  \end{array}\right. \] 
	
	\item  Rectangular pulse:
	
	\[x_{i} (t)=\left\{\begin{array}{lll} 
		(bP/a)(e^{a_{i} t} -1) & if & t\le t_{d}  \\
		x_{i} (t_{d} )e^{a_{i} (t-t_{d} )} & if & t\ge t_{d}  \end{array}\right. \] 
\end{itemize}

\subsection{ Results and discussion}
The figures presented in this section show the numerical resolution of the system of equations of a dumped track section subjected to a rectangular and sinusoidal loads. The properties of the track are defined in table 1.
In figure 6 and figure 7, the sinusoidal pulse is presented in red; however, the rectangular pulse is presented in black. 

\begin{enumerate}
	\item  Displacements and rotations of the rail
	\begin{figure}
		\begin{center}	
			\includegraphics[width=6.3in]{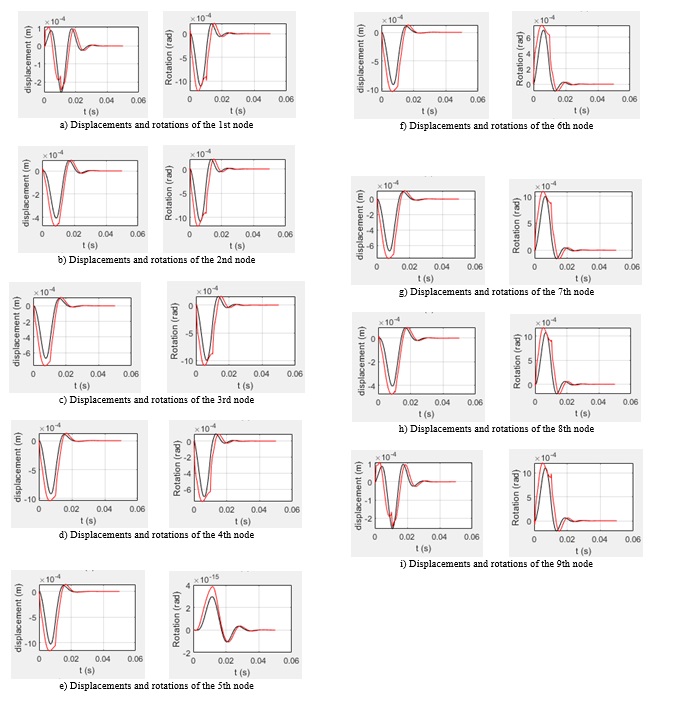} 	
		\end{center}
		\caption{Rail response under sinusoidal, rectangular pulses ($N = 4$, $t_{d} = 0.01$s, $P = 10$T)}	
	\end{figure}
	
	\item 	Displacements of the sleepers
	\begin{figure}
		\begin{center}
			\includegraphics[width=3.8in]{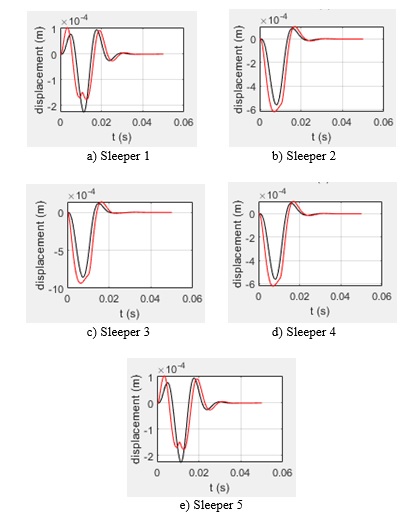}
		\end{center}
		\caption{Sleeper response under sinusoidal and rectangular pulses ($N = 4$, $t_{d} = 0.01$s, $P = 10$T)}
	\end{figure}
\end{enumerate}

\begin{figure}
	\begin{center}
		\includegraphics[scale=0.5]{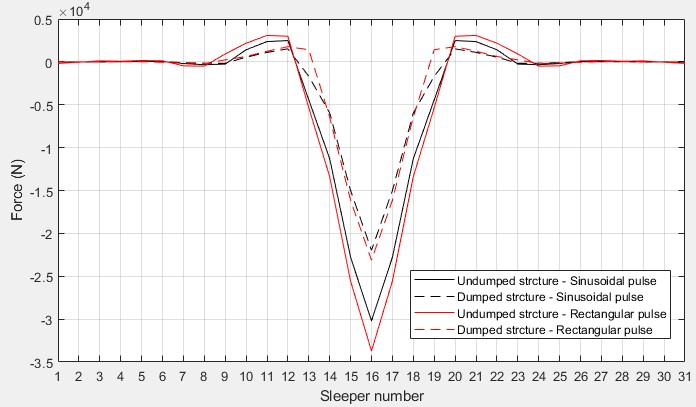}
	\end{center}
	\caption{Loads induced in the substructure ($N = 30$, $t_d = 0.01$s, $P = 10$T)}
\end{figure}

\begin{table}
	\begin{center}
     \begin{tabular}{|p{0.9in}|p{1.1in}|p{1.1in}|p{1.1in}|p{1.1in}|} 
			\hline 	
		Sleeper
		
		 number &   \% of load 
		 
		 (undumped - 
		 
		 sinusoidal load) & \% of load
		 
		  (dumped - 
		  
		  sinusoidal load) & \% of load 
		  
		  (undumped - 
		  
		  rectangular 
		  
		  load) & \% of load 
		  
		  (dumped - 
		  
		  rectangular
		  
		   load) \\
			 \hline  13 & 4.49\% & 1.73\% & 5.50\% & - \\ \hline  \textit{14} & 11.52\% & 6.12\% & 13.66\% & 6.63\% \\
			  \hline \textit{15} & 23.24\% & 15.29\% & 26.00\% & 16.41\% \\ 
			  \hline \textit{16} & 30.78\% & 22.43\% & 34.35\% & 23.55\% \\ 
			  \hline \textit{17} & 23.24\% & 15.29\% & 26.00\% & 16.41\% \\ 
			  \hline  \textit{18} & 11.52\% & 6.12\% & 13.66\% & 6.63\% \\ 
			  \hline \textit{19} & 4.49\% & 1.73\% & 5.50\% & - \\ 
			  \hline 
		\end{tabular}
		\caption{Repartition of the loads under the sleepers ($N = 30$, $t_d = 0.01$s, $P = 10$T)}
	\end{center}
\end{table}

It's clear that the maximum values of rail and sleepers movement under rectangular pulse are higher than those reached under a sinusoidal pulse. The figure 8 shows the maximum loads induced in the substructure. The table 2 shows the repartition of the loads under the sleepers.

These results have many consequences in the railway field. Actually, we may optimize railway infrastructure components for example (like ballast height). Moreover, the study is made by considering a static load (10 T). This load is mainly amplified by rail/wheel interaction and train speed \cite{Key1}.

\section{Conclusion}
Based on the results of the model analysis studied in order to determine the loads induced on the substructure, the following conclusions can be drawn:
\begin{itemize}
	\item 	The common modelling of the load applied on the track due to a moving wheel as a rectangular pulse acting in the time sample of a force signal generates a higher rate of movement in the track and over sizes the loads induced in the substructure than a sinusoidal pulse model;
	\item  Dumping matrix has a major influence on reducing the loads induced in the substructure. Therefore, it’s necessary to preserve the quality of the track components while maintaining it.
\end{itemize}
As an application, we may evaluate the track behavior according to different characteristics of the track elements that degrade because of maintenance operations. Indeed, the ballast is considered as the most affected element because of operations of damping required for track geometry corrections.


\EditInfo{February 22, 2021}{April 10, 2021}{Giuseppe Gaeta}

\end{paper}